\documentclass{article}
\usepackage{amssymb}
\usepackage{amsmath}
\usepackage{abstract}
\usepackage[utf8]{inputenc}
\usepackage[hidelinks]{hyperref}

\usepackage{titlesec}
\titleformat{\section}[block]{\center\scshape\large}{\thesection.}{0.7em}{}
\titleformat{\subsection}[block]{\center}{\thesubsection.}{0.7em}{}
\renewcommand\thesection{\Roman{section}}

\usepackage{amsthm}
\newtheorem{theorem}{Theorem}
\newtheorem{lemma}{Lemma}
\newtheorem{corollary}{Corollary}
\theoremstyle{definition}
\newtheorem*{remark}{Remark}

\newcommand{\keywords}[1]
{
	\small
	\textbf{Keywords:} #1
}

\title{On Goldfeld's Proof of Siegel's Theorem}

\usepackage{authblk}
\author[a]{Zihao Liu}
\affil[a]{International Department, The Affiliated High School of SCNU,\authorcr Email: \url{mailto:travor_lzh@163.com}}

\date{}

\begin{document}
	\maketitle
	\begin{abstract}
		In this paper, we give a detailed account of Goldfeld's proof of Siegel's theorem. Particularly, we present complete proofs of the nontrivial assumptions made in his paper.
	\end{abstract}
	\section{Introduction}
	As stated in \autoref{thsiegel}, Siegel's theorem is a result in multiplicative number theory concerning the lower bound of Dirichlet L-functions associated with quadratic primitive characters.
	\begin{theorem}[Siegel]
		\label{thsiegel}
		For all $\varepsilon>0$ there exists $C(\varepsilon)>0$ such that for any primitive quadratic character $\chi$ modulo $q>1$,
		\begin{equation}
			L(1,\chi)>C(\varepsilon)q^{-\varepsilon}
		\end{equation}
	\end{theorem}
	By applying the mean value theorem for differentiable real-valued functions\footnote{Details can be found in §11.2 of \cite{montgomery_multiplicative_2007} and §21 of \cite{davenport_multiplicative_1980}}, \autoref{thsiegel} allows us to establish a zero-free region\footnote{In fact, every lower bound of $L(1,\chi)$ associated with quadratic $\chi$ can be converted to a certain zero-free region using this argument.} for $L(s,\chi)$:
	\begin{theorem}[Siegel]
		\label{thsiegelz}
		For all $\varepsilon>0$ there exists $C(\varepsilon)>0$ such that for any primitive quadratic character $\chi$ modulo $q>1$, its associated $L(s,\chi)$ is free of zeros in
		\begin{equation}
			s>1-C(\varepsilon)q^{-\varepsilon}
		\end{equation}
	\end{theorem}
	Siegel's theorem in the form of \autoref{thsiegelz} is significant because it provides the prime number theorem in arithmetic progressions with an error term independent of the choice of modulus:
	\begin{theorem}[Siegel-Walfisz\footnote{See §22 of \cite{davenport_multiplicative_1980}}]
		Let $A$ be any fixed positive number and $\pi(x;q,a)$ denote the number of primes that are $\equiv a\pmod q$ and $\le x$. Then for all $(q,a)=1$ we have
		\begin{equation}
			\pi(x;q,a)={1\over\varphi(q)}\int_2^x{\mathrm du\over\log u}+\mathcal O_A\left(x\over\log^Ax\right)
		\end{equation}
	\end{theorem}
	From a historical perspective, \autoref{thsiegel} is an outcome of Siegel's \cite{Siegel1935} investigation on the algebraic properties of the Dedekind zeta function associated with the quartic number field $K=\mathbb Q(\sqrt{q_1},\sqrt{q_2})$:
	\begin{equation}
		\label{eqnfdef}
		f(s)=\zeta(s)L(s,\chi_1)L(s,\chi_2)L(s,\chi_1\chi_2)
	\end{equation}
	where $\chi_1,\chi_2$ denote primitive quadratic characters modulo $q_1,q_2$ respectively.

	Siegel's original proof of \autoref{thsiegel} uses algebraic number theory, but purely analytic proofs have been developed by Estermann\cite{Estermann1949} and Goldfeld\cite{goldfeld_simple_1974}.

	Although Goldfeld's method can lead to a proof of \autoref{thsiegel}, he did not justify all the steps he took in the derivation. Notably, Goldfeld states without justification that
	\begin{equation}
		\label{eqng}
		1\ll{1\over2\pi i}\int_{2-i\infty}^{2+i\infty}f(s+\beta){x^s\over s(s+1)(s+2)(s+3)(s+4)}\mathrm ds
	\end{equation}
	holds for all $0<\beta<1$. The right hand side of \eqref{eqng} is undenably positive, but if the $\ll$ constant relies on the choice of $q_2$ then the subsequent steps in his paper will not lead to a valid proof of \autoref{thsiegel}, so in this paper, we present a more complete version of Goldfeld's proof that addresses these unproven assumptions.

	\section{Plan for the Proof}

	Following the convention, we let $\lambda$ denote the residue of $f(s)$ at $s=1$.
	\begin{equation}
		\lambda=L(1,\chi_1)L(1,\chi_2)L(1,\chi_1\chi_2)
	\end{equation}
	Basically, we attack the problem by giving lower estimate for $\lambda$. When its lower bound is combined with the upper bounds of various $L(1,\chi)$, \autoref{thsiegel} comes out.

	To obtain a positive underestimate of $\lambda$, Goldfeld applied a variant of Perron's formula to $f(s)$, so the lower bound emerges from an application of residue theorem.

	Since Perron's formula connects partial sums with Dirichlet seris, it would be necessary to investigate both of them before proving \autoref{thsiegel}. In \autoref{sfar}, we perform an extensive study of the partial sum associated with $f(s)$ to deduce
	\begin{theorem}
		\label{thmaf}
		Let $a_n$ be the Dirichlet series coefficient of $f(s)$ and $A_0(x,w)$ denote the partial sum of $f(s)$:
		\begin{equation}
			A_0(x,w)=\sum_{n\le x}{a_n\over n^w}
		\end{equation}
		and
		\begin{equation}
			A_k(x,w)=\int_0^xA_{k-1}(y,w)\mathrm dy
		\end{equation}
		Then for every nonnegative integer $k$, there exists an absolute constant $c_k$ and $x_k$ such that $A_k(x,w)\ge c_kx^k$ whenenver $x\ge x_k$.
	\end{theorem}
	Then, in \autoref{sfan}, we study the analytic properties of $f(s)$ to derive a sharp upper bound for $f(s)$ on the right half plane:
	\begin{theorem}
		\label{thmfb}
		For every $\varepsilon>0$ and $\sigma\ge-\varepsilon$. As $|t|\to\infty$, there is
		\begin{equation}
			f(s)\ll_\varepsilon(q_1q_2)^{1+\varepsilon}|t|^{2+\varepsilon}
		\end{equation}
	\end{theorem}
	Finally, in \autoref{sproof}, we combine \autoref{thmaf} and \autoref{thmfb} via Perron's formula to deduce \autoref{thsiegel}.
	\section{Arithmetical properties of $f(s)$}\label{sfar}
	Before investigating the partial sum, we first focus on the properties of $a_n$.
	\begin{lemma}
		$a_1=1$ and $a_n\ge0$ for all $n\in\mathbb Z^+$.
	\end{lemma}
	\begin{proof}
		It follows from the properties of Dirichlet series and Dirichlet convolution that $a_n$ is multiplicative, meaning $a_1=1$. Taking logarithms on both sides, we see
		\begin{equation}
			\log f(s)=\sum_p\sum_{m\ge1}{1\over mp^{ms}}[1+\chi_1(p^m)][1+\chi_2(p^m)]
		\end{equation}
		This indicates that the Dirichlet series coefficients for $\log f(s)$ is nonnegative, so are the those of $f(s)$.
	\end{proof}
	\begin{corollary}
		\label{cola0}
		$A_0(x,w)\ge1$ whenever $x\ge1$.
	\end{corollary}
	Performing an induction on \autoref{cola0}, we can prove \autoref{thmaf}.
	\renewcommand{\proofname}{Proof of \autoref{thmaf}}
	\begin{proof}
		Suppose \autoref{thmaf} is true for $k=m-1\ge0$, then by definition
		\begin{equation}
			\label{eqnam}
			A_m(x,w)\ge c_{m-1}\int_{x_{m-1}}^x y^{m-1}\mathrm dy={c_{m-1}\over m}(x^m-x_{m-1}^m)
		\end{equation}
		The rightmost quantity is $\gg_mx^m$, so we conclude there exists admissible $c_m>0$ and $x_m>0$ such that $A_m(x,w)\ge c_mx^m$ whenever $x\ge x_m$.
	\end{proof}
	\renewcommand{\proofname}{Proof}
	\begin{remark}
		Since constants appearing on the right hand side of \eqref{eqnam} only depends on $m$, \autoref{thmaf} virtually provides positive lower bounds for $A_k(x,w)$ that are independent of $\chi_1$ and $\chi_2$.
	\end{remark}
	\section{Analytic properties of $f(s)$}\label{sfan}
	To obtain an upper estimate for $f(s)$, we quote a classical result from literature:
	\begin{lemma}[Corollary 10.10 of \cite{montgomery_multiplicative_2007}]
		\label{dlf}
		Let $\chi$ be a primitive character modulo $q>1$, and suppose $\sigma$ lies in a fixed interval and $|t|\to\infty$, then
		\begin{equation}
			|L(s,\chi)|\asymp(q|t|)^{1/2-\sigma}|L(1-s,\overline\chi)|
			\footnote{We use the convention $s=\sigma+it$ (where $\sigma,t\in\mathbb R$) throughout the entire paper}
		\end{equation}
	\end{lemma}
	This allows us to conclude that
	\begin{lemma}
		\label{dlb}
		Let $\chi$ be a primitive character modulo $q>1$. For every $\varepsilon>0$ and $\sigma\ge-\varepsilon$, as $|t|\to\infty$ we have
		\begin{equation}
			\label{eqndlbound}
			L(s,\chi)\ll_\varepsilon(q|t|)^{1/2+\varepsilon}
		\end{equation}
	\end{lemma}
	\begin{proof}
		By definition, we see that when $\sigma=-\varepsilon<0$, the Dirichlet series expansion for $L(1-s,\overline\chi)$ converges absolutely. Combining this fact with \autoref{dlf}, we see that \eqref{eqndlbound} holds for $\sigma=-\varepsilon$. Now, it follows from Phragm\'en-Lindel\"of theorem\footnote{See §5.65 of \cite{titchmarsh_theory_2002}} that \eqref{eqndlbound} is valid throughout $\sigma\ge-\varepsilon$ uniformly.
	\end{proof}
	\begin{remark}
		We can develop an argument analogous to the proof of \autoref{dlb} to deduce for all $\sigma\ge0$ and $\varepsilon>0$,
		\begin{equation}
			\label{eqnzb}
			\zeta(s)\ll_\varepsilon|t|^{1/2+\varepsilon}
		\end{equation}
	\end{remark}
	Plugging \autoref{dlb} and \eqref{eqnzb} into \eqref{eqnfdef}, we obtain \autoref{thmaf}.
	\section{Proof of \autoref{thsiegel}}\label{sproof}
	Applying Perron's formula\cite{titchmarsh_theory_1986}, we see that for all $x\ge1$,
	\begin{equation}
		\label{eqna1f}
		A_1(x,w)={1\over2\pi i}\int_{2-i\infty}^{2+i\infty}f(s+w){x^{s+1}\over s(s+1)}\mathrm ds
	\end{equation}
	Integrating on both side of \eqref{eqna1f} for three times and applying \autoref{thmaf}, we have the following result:
	\begin{lemma}[Justification of \eqref{eqng}]
		\label{lmp}
		There exists absolute constants $X,M>0$ such that the following inequality
		\begin{equation}
			M\le{1\over2\pi i}\int_{2-i\infty}^{2+i\infty}f(s+w){x^s\over s(s+1)(s+2)(s+3)(s+4)}\mathrm ds:=J(x,w)
		\end{equation}
		holds uniformly for $x\ge X$.
	\end{lemma}
	To study the integral $J(x,w)$, we set $\lambda=L(1,\chi_1)L(1,\chi_2)L(1,\chi_1\chi_2)$ so that shifting the line of integration to $0>\sigma=-w>-1$ gives
	\begin{align}
		J(x,w)
		&={\lambda x^{1-w}\over(1-w)w(w+1)(w+2)(w+3)} \\
		&+{f(w)\over 4!}\label{eqnfw} \\
		&+\int_{-w-i\infty}^{-w+i\infty}f(s+w){x^s\over s(s+1)(s+2)(s+3)(s+4)}\mathrm ds
	\end{align}
	The upper bound for the remaining integral can be deduced using \autoref{thmfb}:
	\begin{equation}
		\label{eqniw}
		\int_{-w-i\infty}^{-w+i\infty}\ll_\varepsilon{(q_1q_2)^{1+\varepsilon}x^{-w}\over w(1-w)}
	\end{equation}
	Choosing $w$ subtly allows us to omit the contribution of \eqref{eqnfw} when estimating $J(x,w)$:
	\begin{lemma}
		\label{lmin}
		For every $\varepsilon>0$ there exists a primitive quadratic $\chi_1$ modulo $q_1$ and $1-\varepsilon<\beta<1$ such that $f(\beta)\le0$ for all quadratic primitive character $\chi_2$.
	\end{lemma}
	\begin{proof}
		If no quadratic primitive $\chi$ can be found such that $L(s,\chi)$ has a real zero in $(1-\varepsilon,1)$. Then we can use the fact that $\lambda>0$\footnote{See Theorem 4.9 of \cite{montgomery_multiplicative_2007}} and the fact that $\zeta(\sigma)<0$ for $0<\sigma<1$\footnote{See Corollary 1.14 of \cite{montgomery_multiplicative_2007}} to conclude $f(\beta)<0$ for any $1-\varepsilon<\beta<1$.

		If such quadratic primitive $\chi$ does exist, then we let $\chi_1=\chi$ and $\beta$ be the real zero of $L(s,\chi_1)$ in $(1-\varepsilon,1)$ so that $f(\beta)=0$ independent of what $\chi_2$ is. 
	\end{proof}
	\begin{remark}
		\autoref{lmin} explains why the implied constant in \autoref{thsiegel} is not effectively computable because $\chi_1$ and $\beta$ cannot be determined within finitely many steps.
	\end{remark}
	Since $\beta$ and $q_1$ are only associated with $\varepsilon>0$, we can simplify $J(x,w)$ significantly using \autoref{lmin} and \eqref{eqniw} when $w=\beta$:
	\begin{equation}
		\label{eqnjxb}
		J(x,\beta)\le{\lambda x^{1-\beta}\over(1-\beta)(1-\varepsilon)}+\mathcal O_\varepsilon\left(q_2^{1+\varepsilon}x^{-\beta}\right)
	\end{equation}
	\renewcommand{\proofname}{Proof of \autoref{thsiegel}}
	Finally, we can start proving Siegel's theorem. During the proof, $b_1(\varepsilon),b_2(\varepsilon),\dots$ always denote positive constants that only depend on $\varepsilon>0$.
	\begin{proof}
		Plugging \eqref{eqnjxb} into \autoref{lmp}, we know that
		\begin{equation}
			\label{eqnb1}
			b_1(\varepsilon)<\lambda x^{1-\beta}+q_2^{1+\varepsilon}x^{-\beta}
		\end{equation}
		Now we choose $x$ large enough so that the latter term is less than $b_1(\varepsilon)$. This means that we can choose $x$ large enough so that $b_2(\varepsilon)-q_2^{1+\varepsilon}x^{-\beta}>0$, meaning that we can pick
		\begin{equation}
			\label{eqnxb}
			x^\beta=b_3(\varepsilon)q_2^{1+\varepsilon}
		\end{equation}
		where $b_3(\varepsilon)>0$ is a large constant depending on $\varepsilon$. Without loss of generality, we assume $\varepsilon<1/3$, so that plugging \eqref{eqnxb} into \eqref{eqnb1} gives
		\begin{align}
			\lambda
			&>b_2(\varepsilon)x^{-(1-\beta)}=b_4(\varepsilon)q_2^{-(1+\varepsilon)(1-\beta)/\beta} \\
			&> b_4(\varepsilon)q_2^{-(1+\varepsilon)\varepsilon/(1-\varepsilon)}>b_4(\varepsilon)q_2^{-2\varepsilon}
		\end{align}
		To tranfer the lower bound for $\lambda$ to $L(1,\chi_2)$, it suffices to note that for every nontrivial character $\chi$ modulo $q$
		\begin{equation}
			L(1,\chi)=\sum_{n\le T}{\chi(n)\over n}+\mathcal O\left(\frac qT\right)
		\end{equation}
		as it indicates $L(1,\chi)\ll\log q$ after setting $T=q$. Consequently, we have
		\begin{align}
			L(1,\chi_2)
			&>b_5(\varepsilon)q_2^{-2\varepsilon}(\log q_1)^{-1}(\log q_1q_2)^{-1} \\
			&>b_6(\varepsilon)q_2^{-2\varepsilon}(\log q_1q_2)^{-1}\label{eqnps}
		\end{align}
		If we make $q_2$ be sufficiently large, then \eqref{eqnps} gets simplified into
		\begin{equation}
			\label{eqnl1x}
			L(1,\chi_2)>b_7(\varepsilon)q_2^{-2\varepsilon}(\log q_2)^{-1}>b_7(\varepsilon)q_2^{-3\varepsilon}
		\end{equation}
		Finally, we make $b_7(\varepsilon)$ very small to ensure that \eqref{eqnl1x} hold for small values of $q_2$, so the proof of \autoref{thsiegel} is complete.
	\end{proof}
	\bibliographystyle{plain}
	\bibliography{refs.bib}

\begin{thebibliography}{1}

\bibitem{davenport_multiplicative_1980}
Harold Davenport.
\newblock {\em Multiplicative {Number} {Theory}}, volume~74 of {\em Graduate
  {Texts} in {Mathematics}}.
\newblock Springer New York, New York, NY, 1980.

\bibitem{Estermann1949}
T.~{Estermann}.
\newblock {On Dirichlet's \(L\) functions}.
\newblock {\em {J. Lond. Math. Soc.}}, 23:275--279, 1949.

\bibitem{goldfeld_simple_1974}
D.~M. Goldfeld.
\newblock A {Simple} {Proof} of {Siegel}'s {Theorem}.
\newblock {\em Proceedings of the National Academy of Sciences},
  71(4):1055--1055, April 1974.

\bibitem{montgomery_multiplicative_2007}
Hugh~L. Montgomery and Robert~C. Vaughan.
\newblock {\em Multiplicative number theory {I}: classical theory}.
\newblock Number~97 in Cambridge studies in advanced mathematics. Cambridge
  University Press, Cambridge, UK ; New York, 2007.
\newblock OCLC: ocm61757122.

\bibitem{Siegel1935}
Carl Siegel.
\newblock Über die classenzahl quadratischer zahlkörper.
\newblock {\em Acta Arithmetica}, 1(1):83--86, 1935.

\bibitem{titchmarsh_theory_2002}
E.~C. Titchmarsh.
\newblock {\em The theory of functions}.
\newblock Oxford science publications. Oxford Univ. Press, Oxford, 2. ed.,
  reprinted edition, 2002.
\newblock OCLC: 249703508.

\bibitem{titchmarsh_theory_1986}
E.~C. Titchmarsh and D.~R. Heath-Brown.
\newblock {\em The theory of the {Riemann} zeta-function}.
\newblock Oxford science publications. Oxford University Press, New York, 2nd
  ed edition, 1986.

\end{thebibliography}

\end{document}